\documentclass{amsart}

\usepackage{amssymb}
\usepackage{amsthm}

\usepackage{amsfonts}
\usepackage{mathrsfs}
\usepackage{amsmath}

\newtheorem{theorem}{Theorem}
\newtheorem{lemma}[theorem]{Lemma}
\newtheorem{proposition}[theorem]{Proposition}
\newtheorem{thm}{Theorem}
\newtheorem*{conj}{Conjecture}

\newtheorem{remark}[thm]{Remark}

\newtheorem*{thma}{Theorem A}
\newtheorem*{thmb}{Theorem B}

\newcommand{\biggg}{\bBigg@{3}}

\newcommand{\Biggg}{\bBigg@{3.5}}

\newcommand{\bigggg}{\bBigg@{4}}

\newcommand{\Bigggg}{\bBigg@{4.5}}

\newcommand{\ball}{\mathbb{B}_n}

\newcommand{\bbR}{\mathbb{R}}

\newcommand{\hyperg}[4]{\: _2\! F_1\! \left[\!\! \begin{array}{c} #1,\, #2 \\[4pt] #3 \end{array}\!\! ;\,  #4 \right]}

\title{Schwarz-Pick lemma for harmonic functions}

\thanks{This work was supported by the National Natural Science
Foundation of China grant 11971453.}

\author{Congwen Liu}
\email{cwliu@ustc.edu.cn}

\address{School of Mathematical Sciences,
University of Science and Technology of China\\
Hefei, Anhui 230026,
People's Republic of China\\
and\\
CAS Wu Wen-Tsun Key Laboratory of Mathematics\\
USTC, Hefei, China}


\begin{document}

\begin{abstract}
Based on the recently proved Khavinson conjecture, we establish an inequality of Schwarz-Pick type for
harmonic functions on the unit ball of $\mathbb{R}^n$.
\end{abstract}

\maketitle

\section{Introduction}

This is a sequel to the paper \cite{Liu20}. There we proved the Khavinson conjecture,
which says for bounded harmonic functions on the unit ball of $\bbR^n$, the sharp constants in
the estimates for their radial derivatives and for their gradients coincide.
In this paper, we further prove the following
\begin{theorem}\label{thm:main}
Let $u$ be a real-valued harmonic function on the unit ball $\ball$ of $\mathbb{R}^n$ and
$|u| < 1$ on $\ball$.
\begin{enumerate}
\item[(i)]
When $n=2$ or $n\geq 4$, the following sharp inequality holds:
\begin{equation}\label{eqn:main}
|\nabla u(x)| \leq \frac {2m_{n-1}(\mathbb{B}_{n-1})}{m_n(\ball)} \frac{1}{1-|x|^2}, \qquad x\in \ball,
\end{equation}
where $m_n$ denotes the Lebesgue measure on $\mathbb{R}^n$.
Equality holds if and only if $x=0$ and $u = U \circ T$ for some orthogonal transformation $T$,
where $U$ is the Poisson integral of the function that equals $1$ on a hemisphere and $-1$ on the remaining
hemisphere.
\item[(ii)]
When $n=3$ we have
\begin{equation}\label{eqn:maincase3}
|\nabla u(x)| < \frac {8}{3\sqrt{3}} \frac{1}{1-|x|^2}, \qquad  x\in \mathbb{B}_3.
\end{equation}
The constant $\frac {8}{3\sqrt{3}}$ here is the best possible.
\end{enumerate}
\end{theorem}

\begin{remark}
Curiously, the inequality \eqref{eqn:main} fails when $n=3$. Note that $\frac {8}{3\sqrt{3}}\approx 1.5396$, while
the constant $\frac {2m_{n-1}(\mathbb{B}_{n-1})}{m_n(\ball)}$ in \eqref{eqn:main} equals $\frac {3}{2}$ when $n=3$.
According to \cite{Kha92}, given $x_0\in \mathbb{B}_3$, there is a $u_0$ harmonic in $\mathbb{B}_3$,
$|u_0|<1$, satisfying
\[
\left|\left\langle\nabla u_0(x_0), \frac {x_0}{|x_0|}\right\rangle\right|
~=~ \frac {(9-|x_0|^2)^2} {3\sqrt{3} (1-|x_0|^2) [(|x_0|^2+3)^{3/2}+ 3\sqrt{3} (1-|x_0|^2)]}.
\]
Since
\[
\lim_{t\to 1-} \frac {(9-t^2)^2} {3\sqrt{3} [(t^2+3)^{3/2}+ 3\sqrt{3} (1-t^2)]} = \frac {8}{3\sqrt{3}} > \frac {3}{2},
\]
we see that if $x_0$ is sufficiently near the boundary $\mathbb{S}^2$ then $|\nabla u_0(x_0)|> \frac {3}{2} \frac {1}{1-|x_0|^2}$, 
which contradicts the inequality \eqref{eqn:main}. See Remark 2 in the next section for more explanations.
\end{remark}

We refer to Theorem \ref{thm:main} as the Schwarz-Pick lemma for harmonic functions, for
(1) it is in analogy to a weaker version of the classical Schwarz-Pick lemma, which states that every
holomorphic mapping $f$ of the unit disk onto itself satisfies the inequality:
\[
|f^{\prime}(z)| \leq \frac {1}{1-|z|^2}, \qquad |z|<1;
\]
(2) it is a generalization of the following harmonic Schwarz lemma:
\begin{thma}[{\cite[Theorem 6.26]{ABR01}}]
If $u$ is a real-valued harmonic function on $\ball$ and $|u| < 1$ on $\ball$, then
\begin{equation}
\label{eqn:ABR01}
|\nabla u(0)| \leq \frac {2m_{n-1}(\mathbb{B}_{n-1})}{m_n(\ball)}.
\end{equation}
Equality holds if and only if $u = U \circ T$ for some orthogonal transformation $T$,
where $U$ is as in Theorem \ref{thm:main}.
\end{thma}

Some special cases of Theorem \ref{thm:main} are known.
When $n=2$, the inequality \eqref{eqn:main} reads
\begin{equation}\label{eqn:Lindelof}
|\nabla u(z)| \leq  \frac{4}{\pi }\frac{1}{1-|z|^2}, \qquad |z|<1,
\end{equation}
which is a reformulation of the Lindel\"of inequality in the unit disc (see \cite[Theorem 3]{Col89}).
Recently, Kalaj \cite{Kal17} established the inequality \eqref{eqn:main} for $n=4$.
When $n=3$, although not explicitly stated, the inequality \eqref{eqn:maincase3} is an easy consequence of the main result of
\cite{Mel19}.

When $n=2$ or $n\geq 4$, Theorem \ref{thm:main} can be restated as
\begin{equation*}\label{eqn:main2}
|\nabla u(x)| \leq \frac {2m_{n-1}(\mathbb{B}_{n-1})}{m_n(\ball)} \frac{1}{1-|x|^2}
\sup_{y\in \ball} |u(y)|
\end{equation*}
for bounded harmonic functions $u$ on $\ball$.
This is obviously related to the following classical estimate (see for instance \cite[p.139, (6)]{PW84})
\begin{equation}
\label{eqn:PW6}
|\nabla u(x)|\leq \frac {m_{n-1}(\mathbb{B}_{n-1})}{m_n(\ball)} \frac{1}{d(x)} \operatorname{osc}_{\Omega}(u)
\end{equation}
for harmonic functions $u$ in $\Omega\subset \mathbb{R}^n$, where $\operatorname{osc}_{\Omega}(u)$ is the
oscillation of $u$ in $\Omega$ and $d(x)$ denotes the distance of $x \in
\Omega$ to the boundary $\partial\Omega$.
In particular, if $\Omega=\ball$ and $u$ is a harmonic function on $\ball$ with $|u| < 1$, then \eqref{eqn:PW6} reads
\begin{equation*} 
|\nabla u(x)| \leq \frac {2m_{n-1}(\mathbb{B}_{n-1})}{m_n(\ball)} \frac{1}{1-|x|}.
\end{equation*}
Compare this with the inequality \eqref{eqn:main}.

It is also interesting to compare the inequality \eqref{eqn:main2} with the following sharp inequality
in \cite{KM10b} (see also \cite[p.131]{KM12}):
\begin{equation}\label{eqn:KMhalfspace}
|\nabla v (x)|\leq
\frac {4\, (n-1)^{\frac {n+1}{2}}\, m_{n-1}(\mathbb{B}_{n-1})} { n^{\frac {n+2}{2}}\, m_n(\ball)}
\frac {1}{x_n} \sup_{y\in \mathbb{R}_{+}^n} |v(y)|.
\end{equation}
Here, $v$ is a bounded harmonic function in the half--space $\mathbb{R}^n_+:=\{(x^{\prime},x_n)\in \mathbb{R}^{n}:x_n>0\}$.

Recently,  several versions of Schwarz lemma for harmonic functions or harmonic mappings
were established. See \cite{Che13, Kal16, KV12, Mar15, Mat18, Mel18}. 
In particular, Kalaj and Vuorinen \cite{KV12} obtained the following refinement of 
the inequality \eqref{eqn:Lindelof}:
\begin{equation*}\label{eqn:KV12}
|\nabla u(z)| \leq~  \frac{4}{\pi }\frac{1-|u(z)|^2}{1-|z|^2}, \qquad |z|<1.
\end{equation*}
This, together with our Theorem \ref{thm:main},  suggests the following
\begin{conj}
Under the hypotheses of Theorem \ref{thm:main}(i), we have
\begin{equation*}
|\nabla u(x)| \leq \frac {2m_{n-1}(\mathbb{B}_{n-1})}{m_n(\ball)} \frac{1-|u(x)|^2}{1-|x|^2}, \qquad x\in \ball.
\end{equation*}
\end{conj}
We are not able to prove this conjecture and leave it as an open question.

\section{Outline of the proof}

Since the cases $n=2$ and $n=3$ are known, we shall prove only Theorem \ref{thm:main} for $n\geq 4$.

Recall that we proved in \cite{Liu20} the following
\begin{thmb}[{\cite[Theorem 2]{Liu20}}] 
Let $n\geq 3$ and let $u$ be a real-valued bounded harmonic function on $\ball$. We have the
following sharp inequality:
\[
|\nabla u(x)| \leq C(x) \sup_{y\in \ball} |u(y)|, \qquad x\in \ball,
\]
with
\[
C(x):=\frac {(n-1) m_{n-1}(\mathbb{B}_{n-1})}{m_n(\ball)} \left\{\int\limits_{-1}^1
\frac {\left| t - \frac{n-2}n |x| \right| (1-t^2)^{\frac {n-3}{2}}}
{(1-2t|x|+|x|^2)^{\frac {n-2}{2}}} dt \right\} \frac{1}{1-|x|^2}.
\]
\end{thmb}
Thus, in order to prove the inequality \eqref{eqn:main}, it suffices to prove the following

\begin{proposition}\label{pro:monotonicity}
When $n\geq 4$, the function
\begin{equation*}\label{eqn:funPhi}
\Phi(\rho) := \int\limits_{-1}^1 \frac {\left| t - \frac{n-2}n \rho \right| (1-t^2)^{\frac {n-3}{2}}}
{(1 - 2 t\rho + \rho^2)^{\frac {n-2}{2}}} dt
\end{equation*}
is strictly decreasing on $[0, 1]$ and
\[
\max_{\rho\in [0,1]} \Phi(\rho) = \Phi(0) =\frac {2}{n-1}.
\]
\end{proposition}

\begin{remark}
In contrast to the case $n\geq 4$, if $n=3$ then
\[
\Phi(\rho)
= \frac {2}{3} \frac {(1+\frac {1}{3}\rho^2)^{\frac {3}{2}}-1+\rho^2}{\rho^2}
\]
is strictly increasing on $[0, 1]$ and attains its maximum at $\rho=1$.
This explains why the inequality \eqref{eqn:main} fails when $n=3$, as well as
why, unlike \eqref{eqn:main}, the inequality \eqref{eqn:maincase3} is sharp but
always strict.
\end{remark}

Assuming Proposition \ref{pro:monotonicity} for the moment, we shall prove the second 
assertion of Theorem \ref{thm:main}(i).
In view of Theorem B, it follows from the strict monotonicity of $\Phi$ that 
the equality in \eqref{eqn:main} takes place if and only if $x=0$. Then, by
Theorem A, $u$ must be of the form $u = U \circ T$, with $T$ an orthogonal transformation
and $U$ the Poisson integral of the function that equals $1$ on a hemisphere and $-1$ on the remaining
hemisphere.

We now turn to the proof of Proposition \ref{pro:monotonicity}.
An easy computation leads to $\Phi^{\prime}(0)=0$. Thus, the problem is further
reduced to the following
\begin{theorem}\label{pro:concave}
If $n\geq 4$ then $\Phi$ is a strictly concave function on the interval $(0, 1)$.
\end{theorem}

We divide the proof of Theorem \ref{pro:concave} into the following two propositions.

\begin{proposition}\label{prop:2ndderofPhi}
We have
\begin{align}\label{eqn:2ndderofPhi}
\Phi^{\prime\prime}(\rho) &~=~  \frac {2(n-2)}{\rho^2} \left[1-\frac {(n-2)^2}{n^2} \rho^2 \right]^{\frac {n-3}{2}}
\left(1-\frac {n-4}{n} \rho^2 \right)^{-\frac {n}{2}} \\
&\qquad \times \Biggl\{ \left[1-\frac {(n-2)(n-3)}{n^2} \rho^2 \right] \left(1-\frac {n-4}{n} \rho^2 \right) \notag\\
& \qquad \qquad \quad - \left[1-\frac {(n-2)(n-3)}{n(n-1)} \rho^2 \right] \left(1-\frac {(n-2)^2}{n^2} \rho^2 \right) \notag\\
& \qquad \qquad \qquad \quad \times \left(1-\frac {n-2}{n} \rho^2 \right) \hyperg {1} {\frac {n}{2}} {\frac {n+1}{2}}
 {\frac {\rho^2 \left[1-\frac {(n-2)^2}{n^2} \rho^2 \right]} {1-\frac {n-4}{n} \rho^2}} \Biggr\}.\notag
\end{align}

\end{proposition}

Here and throughout the paper, $\hyperg{a}{b}{c}{z}$ denotes the Gauss hypergeometric function defined by
\begin{equation*}\label{eq:hypergdefin}
\hyperg{a}{b}{c}{z} := \sum_{k=0}^{\infty}\frac{(a)_k(b)_k}{(c)_k\, k!} z^k
\end{equation*}
for $|z|<1$, with $(\lambda)_k$ the Pochhammer symbol (or the
extended factorial), which is defined by
\[
(\lambda)_0 := 1, \quad (\lambda)_k := \lambda(\lambda+1)\ldots(\lambda+k-1)
\quad \text{ for } k\geq 1.
\]

\begin{proposition}\label{lem:technical}
If $n \geq 4$ then
\begin{equation}\label{eqn:technical}
\hyperg {1} {\frac {n}{2}} {\frac {n+1}{2}} {\frac {t \left[1-\frac {(n-2)^2}{n^2} t\right]} {1-\frac {n-4}{n} t}}
 ~>~ \frac {\left( 1-\frac {n-4}{n} t \right) \left[ 1-\frac {(n-2)(n-3)}{n^2} t \right]}
{\left(1-\frac {n-2}{n} t\right) \left[ 1-\frac {(n-2)^2}{n^2} t \right] \left[ 1-\frac {(n-2)(n-3)}{n(n-1)} t \right]}
\end{equation}
holds for all $t\in [0,1]$. 
\end{proposition}

\begin{remark}
As we will see in the proof, if $n=3$ then the inequality in \eqref{eqn:technical} is reversed.
\end{remark}

Propositions \ref{prop:2ndderofPhi} and \ref{lem:technical} will be proved in the next two sections.

\section{The proof of Proposition \ref{prop:2ndderofPhi}}

The proof will be divided into three steps.
\subsection*{Step 1}
We express the function $\Phi$ as follows.
\begin{align}\label{eqn:reprn}
\Phi(\rho) ~=~& \int\limits_{-1}^1 \left| t- \frac {n-2}{n} \rho  \right| (1-t^2)^{\frac {n-3}{2}} dt \\
&\quad +~  (n-2)\rho  \int\limits_{-1}^1 \left| t- \frac {n-2}{n} \rho \right| (1-t^2)^{\frac {n-3}{2}} t dt \notag\\
& \quad +~ \sum_{k=2}^{\infty} \frac {2n(n-2)} {k(k-1)(k+n-2)(k+n-1)} \notag \\
& \quad \qquad \quad \times \left[1-\frac {(n-2)^2}{n^2} \rho^2 \right]^{\frac {n+1}{2}} C_{k-2}^{\frac {n+2}{2}} \left(\frac {n-2}{n} \rho\right) \rho^k, \notag
\end{align}
where $C_k^{\lambda}(x)$ is the Gegenbauer polynomial (also known as the ultraspherical
polynomials) of degree $k$ associated to $\lambda$, which is defined by the generating relation
\begin{equation}\label{eqn:generatingformula}
(1-2xz+z^2)^{-\lambda} = \sum_{k=0}^{\infty} C_k^{\lambda}(x) z^k, \qquad -1<x<1,\; |z|<1 .
\end{equation}

Using the generating relation \eqref{eqn:generatingformula} and noting that
\[
C_k^{\frac{n-2}{2}} (t) \equiv 1 \quad \text{and} \quad  C_k^{\frac{n-2}{2}} (t) = (n-2) t,
\]
we obtain
\begin{align*}
\Phi(\rho) ~=~& \sum_{k=0}^{\infty} \left\{ \int\limits_{-1}^{1} \left| t - \frac {n-2}{n} \rho \right|
(1-t^2)^{\frac{n-3}2} C_k^{\frac{n-2}{2}} (t)  dt \right\} \rho^k\\
=~& \int\limits_{-1}^1 \left| t- \frac {n-2}{n} \rho  \right| (1-t^2)^{\frac {n-3}{2}} dt \\
& \qquad +  (n-2)\rho  \int\limits_{-1}^1 \left| t- \frac {n-2}{n} \rho \right| (1-t^2)^{\frac {n-3}{2}} t dt \notag\\
& \qquad + \sum_{k=2}^{\infty} \left\{ \int\limits_{-1}^{1} \left| t - \frac {n-2}{n} \rho \right|
(1-t^2)^{\frac{n-3}2} C_k^{\frac{n-2}{2}} (t)  dt \right\} \rho^k.
\end{align*}
Then \eqref{eqn:reprn} follows by an application of Lemma \ref{lem:intbyterms} below, 
with $\lambda=\frac {n-2}{2}$ and $s=\frac {n-2}{n}\rho$.
\begin{lemma}[{\cite[Lemma 5]{Liu20}}] \label{lem:intbyterms}
Let $\lambda>-1/2$ and $-1<s<1$. Then we have
\begin{align}
\label{eqn:intbyterms}
\int\limits_{-1}^{1} & |x-s| (1-x^2)^{\lambda-\frac {1}{2}} C_k^{\lambda}(x) dx \\
&~=~ \frac {8\lambda(\lambda+1)} {k(k-1)(k+2\lambda)(k+2\lambda+1)}
(1-s^2)^{\lambda+\frac {3}{2}}  C_{k-2}^{\lambda+2} (s) \notag
\end{align}
for  $k=2,3,\ldots$.
\end{lemma}

\subsection*{Step 2}
We claim the formula
\begin{align}\label{eqn:Phi2ndder1}
\Phi^{\prime\prime}(\rho)
~=~& \frac {2(n-2)^2}{n^2} \left[1-\frac {(n-2)^2}{n^2} \rho^2 \right]^{\frac {n-3}{2}} \sum_{k=0}^{\infty}
C_k^{\frac {n-2}{2}} \left(\frac {n-2}{n} \rho \right) \, \rho^k \\
&\quad -~ \frac {4 (n-2)^2}{n(n-1)} \left[1-\frac {(n-2)^2}{n^2} \rho^2 \right]^{\frac {n-1}{2}}
\sum_{k=0}^{\infty} \frac {(n-1)_k}{(n)_k} C_k^{\frac {n}{2}} \left(\frac {n-2}{n} \rho \right) \, \rho^k \notag\\
&\quad +~ \frac {2(n-2)}{n+1}  \left[1-\frac {(n-2)^2}{n^2} \rho^2 \right]^{\frac {n+1}{2}}
\sum_{k=0}^{\infty} \frac {(n)_k}{(n+2)_k} C_k^{\frac {n+2}{2}} \left(\frac {n-2}{n} \rho \right) \, \rho^k. \notag
\end{align}

The proof is similar to that of Lemma 8 in \cite{Liu20}. In view of \eqref{eqn:reprn}, we write
\begin{align*}
\Phi_1(\rho) :=& \int\limits_{-1}^1 \left| \frac {n-2}{n} \rho  -x \right| (1-x^2)^{\frac {n-3}{2}} dx, \\
\Phi_2(\rho) :=& (n-2)\rho  \int\limits_{-1}^1 \left| \frac {n-2}{n} \rho  -x \right| (1-x^2)^{\frac {n-3}{2}} x dx, \notag\\
\Phi_3(\rho) :=&  \sum_{k=2}^{\infty} \frac {2n(n-2)}{k(k-1)(k+n-2)(k+n-1)} \notag\\
& \hspace{1cm} \times \left\{\left[1-\frac {(n-2)^2}{n^2} \rho^2 \right]^{\frac {n+1}{2}}
C_{k-2}^{\frac {n+2}{2}} \left(\frac {n-2}{n} \rho \right) \right\} \rho^k.
\end{align*}
Straightforward computations yield
\[
\Phi_1^{\prime\prime}(\rho) = \frac {2(n-2)^2}{n^2} \left[1-\frac {(n-2)^2}{n^2} \rho^2 \right]^{\frac {n-3}{2}}
C_0^{\frac {n-2}{2}} \left(\frac {n-2}{n} \rho\right)
\]
and
\begin{align*}
\Phi_2^{\prime\prime}(\rho)
=~& \frac {2(n-2)^2}{n^2} \left[1-\frac {(n-2)^2}{n^2} \rho^2 \right]^{\frac {n-3}{2}}
C_1^{\frac {n-2}{2}} \left(\frac {n-2}{n} \rho\right) \, \rho \\
& \quad -~ \frac {4(n-2)^2}{n(n-1)} \left[1-\frac {(n-2)^2}{n^2} \rho^2 \right]^{\frac {n-1}{2}}
C_0^{\frac {n}{2}} \left(\frac {n-2}{n} \rho\right).
\end{align*}
So, it remains to show that
\begin{align}\label{eqn:2ndderPhi3}
\Phi_3^{\prime\prime}(\rho)
~=&~ \frac {2(n-2)^2}{n^2} \left[1-\frac {(n-2)^2}{n^2} \rho^2 \right]^{\frac {n-3}{2}} \sum_{k=2}^{\infty}
C_k^{\frac {n-2}{2}} \left(\frac {n-2}{n} \rho \right) \, \rho^k \\
&\quad -~ \frac {4 (n-2)^2}{n(n-1)} \left[1-\frac {(n-2)^2}{n^2} \rho^2 \right]^{\frac {n-1}{2}}
\sum_{k=1}^{\infty} \frac {(n-1)_k}{(n)_k} C_k^{\frac {n}{2}} \left(\frac {n-2}{n} \rho \right) \, \rho^k \notag\\
&\quad +~ \frac {2(n-2)}{n+1}  \left[1-\frac {(n-2)^2}{n^2} \rho^2 \right]^{\frac {n+1}{2}}
\sum_{k=0}^{\infty} \frac {(n)_k}{(n+2)_k} C_k^{\frac {n+2}{2}} \left(\frac {n-2}{n} \rho \right) \, \rho^k. \notag
\end{align}
We need the following lemma.
\begin{lemma}[{\cite[Lemma 4]{Liu20}}] \label{lem:diffofgeng2}
If $\lambda\neq 1$ then
\begin{equation*}\label{eqn:diffofgeng2}
\frac {d}{dx} \left\{(1-x^2)^{\lambda-\frac {1}{2}} C_k^{\lambda}(x) \right\}
~=~  - \frac {(k+1)(k+2\lambda-1)}{2(\lambda-1)} (1-x^2)^{\lambda-\frac {3}{2}} C_{k+1}^{\lambda-1}(x).
\end{equation*}
\end{lemma}
Applying Lemma \ref{lem:diffofgeng2} successively, it follows that
\begin{align*}
\frac {d}{d\rho} & \left\{\left[1-\frac {(n-2)^2}{n^2} \rho^2 \right]^{\frac {n+1}{2}}
C_{k-2}^{\frac {n+2}{2}} \left(\frac {n-2}{n} \rho \right) \right\}\\
&\quad = - \frac {n-2}{n} \frac {(k-1) (k+n-1)} {n} \left[1-\frac {(n-2)^2}{n^2} \rho^2 \right]^{\frac {n-1}{2}}
C_{k-1}^{\frac {n}{2}} \left(\frac {n-2}{n} \rho \right)
\end{align*}
and
\begin{align*}
\frac {d^2}{d\rho^2} & \left\{\left[1-\frac {(n-2)^2}{n^2} \rho^2 \right]^{\frac {n+1}{2}}
C_{k-2}^{\frac {n+2}{2}} \left(\frac {n-2}{n} \rho \right) \right\}\\
&\quad =~ \frac {(n-2)^2}{n^2} \frac {k (k+n-2) (k-1) (k+n-1)} {n (n-2)} \\
& \qquad \qquad \times \left[1-\frac {(n-2)^2}{n^2} \rho^2 \right]^{\frac {n-3}{2}}
C_{k}^{\frac {n-2}{2}} \left(\frac {n-2}{n} \rho \right),
\end{align*}
and hence
\begin{align*} 
\Phi_3^{\prime\prime}(\rho) &~=~  \sum_{k=2}^{\infty} \frac {2n(n-2)}{k(k-1)(k+n-2)(k+n-1)} \notag\\
& \quad \times \Biggl( \frac {d^2}{d\rho^2} \biggl\{\left[1-\frac {(n-2)^2}{n^2} \rho^2 \right]^{\frac {n+1}{2}}
C_{k-2}^{\frac {n+2}{2}} \biggl(\frac {n-2}{n} \rho \biggr) \biggr\}\, \rho^k\\
& \quad \hspace{1cm} +~ 2 \frac {d}{d\rho} \biggl\{\left[1-\frac {(n-2)^2}{n^2} \rho^2 \right]^{\frac {n+1}{2}}
C_{k-2}^{\frac {n+2}{2}} \left(\frac {n-2}{n} \rho \right) \biggr\} \, \frac {d}{d\rho} \left(\rho^k\right)\\
& \quad \hspace{1cm} +~ \left[1-\frac {(n-2)^2}{n^2} \rho^2 \right]^{\frac {n+1}{2}}
C_{k-2}^{\frac {n+2}{2}} \left(\frac {n-2}{n} \rho \right)\, \frac {d^2}{d\rho^2} \left(\rho^k\right) \Biggr)\\
&~=~ \frac {2(n-2)^2}{n^2} \left[1-\frac {(n-2)^2}{n^2} \rho^2 \right]^{\frac {n-3}{2}} \sum_{k=2}^{\infty}
C_k^{\frac {n-2}{2}} \left(\frac {n-2}{n} \rho \right) \, \rho^k \\
&\quad - \frac {4 (n-2)^2}{n} \left[1-\frac {(n-2)^2}{n^2} \rho^2 \right]^{\frac {n-1}{2}} \sum_{k=2}^{\infty} \frac {1}{k+n-2}
C_{k-1}^{\frac {n}{2}} \left(\frac {n-2}{n} \rho \right) \, \rho^{k-1} \notag \\
&\quad + \left[1-\frac {(n-2)^2}{n^2} \rho^2 \right]^{\frac {n+1}{2}} \sum_{k=2}^{\infty} \frac {2n(n-2)}{(k+n-2)(k+n-1)}
C_{k-2}^{\frac {n+2}{2}} \left(\frac {n-2}{n} \rho \right) \, \rho^{k-2}.
\end{align*}
This is precisely the desire identity \eqref{eqn:2ndderPhi3}, in view of that
\[
\frac {(n-1)_k}{(n)_k} = \frac {n-1}{n+k-1} \quad \text{and} \quad \frac {(n)_k}{(n+2)_k} = \frac {n(n+1)}{(n+k)(n+k+1)}.
\]

\subsection*{Step 3}
We are now ready to conclude the proof of the Proposition \ref{prop:2ndderofPhi}.

Denote the three terms in the right hand side of \eqref{eqn:Phi2ndder1} by $I$, $II$ and $III$ respectively.
By the generating relation \eqref{eqn:generatingformula}, we have
\begin{align}
I ~=~& \frac {2(n-2)^2}{n^2} \left[1-\frac {(n-2)^2}{n^2} \rho^2 \right]^{\frac {n-3}{2}}
\left[1- 2 \left(\frac {n-2}{n} \rho\right)\rho + \rho^2 \right]^{-\frac {n-2}{2}} \label{eqn:termI}\\
~=~& \frac {2(n-2)^2}{n^2} \left[1-\frac {(n-2)^2}{n^2} \rho^2 \right]^{\frac {n-3}{2}}
\left(1- \frac {n-4}{n} \rho^2 \right)^{-\frac {n-2}{2}}. \notag
\end{align}
For $II$, we shall make use of the following
\begin{lemma}[{\cite[p.279, (8)]{Rai60}}] Suppose that $\nu, \lambda\in \mathbb{R}$ and $2\lambda \neq 0, -1, -2, \ldots$.
Then the identity
\begin{equation}\label{eqn:generatingformula2}
(1-xz)^{-\nu} \hyperg {\frac {\nu}{2}} {\frac {\nu+1}{2}} {\lambda+\frac {1}{2}} {\frac {z^2(x^2-1)}{(1-xz)^2}}
~=~ \sum_{k=0}^{\infty} \frac {(\nu)_k}{(2\lambda)_k} C_k^{\lambda}(x) z^k
\end{equation}
holds, whenever both sides make sense.
\end{lemma}

We apply \eqref{eqn:generatingformula2} with $\nu=n-1$ and $\lambda=\frac {n}{2}$ to obtain
\begin{align}\label{eqn:termII0}
II ~=~ & - \frac {4 (n-2)^2}{n(n-1)} \left[1-\frac {(n-2)^2}{n^2} \rho^2 \right]^{\frac {n-1}{2}}
\left(1-\frac {n-2}{n} \rho^2 \right)^{-n+1}\\
& \quad \times \hyperg {\frac {n-1}{2}} {\frac {n}{2}} {\frac {n+1}{2}}
 {\frac {\rho^2 \left[\frac {(n-2)^2}{n^2} \rho^2 -1\right]} {(1-\frac {n-2}{n} \rho^2 )^2}}. \notag
\end{align}
On the other hand,
\begin{align}\label{eqn:hyperg1}
&\hyperg {\frac {n-1}{2}} {\frac {n}{2}} {\frac {n+1}{2}}
 {\frac {\rho^2 \left[\frac {(n-2)^2}{n^2} \rho^2 -1\right]} {(1-\frac {n-2}{n} \rho^2 )^2}}\\
&\qquad \hspace{1cm} ~=~ \frac {\left(1-\frac {n-2}{n} \rho^2 \right)^n}
{\left(1-\frac {n-4}{n} \rho^2 \right)^{\frac {n}{2}}}
\hyperg {1} {\frac {n}{2}} {\frac {n+1}{2}}
 {\frac {\rho^2 \left[1-\frac {(n-2)^2}{n^2} \rho^2 \right]} {1-\frac {n-4}{n} \rho^2}}, \notag
\end{align}
which follows from Pfaff's transformation formula (\cite[p.105, (4)]{EMOT53a})
\[
\hyperg{a}{b}{c}{z} = (1-z)^{-b} \hyperg{c-a}{b}{c}{\frac {z}{z-1}}.
\]
Substituting \eqref{eqn:hyperg1} into \eqref{eqn:termII0} yields
\begin{align}\label{eqn:termII}
II ~=~ & - \frac {4 (n-2)^2}{n(n-1)} \left[1-\frac {(n-2)^2}{n^2} \rho^2 \right]^{\frac {n-1}{2}}
\left(1-\frac {n-2}{n} \rho^2 \right) \\
& \quad \times \left(1-\frac {n-4}{n} \rho^2 \right)^{-\frac {n}{2}}
\hyperg {1} {\frac {n}{2}} {\frac {n+1}{2}}
 {\frac {\rho^2 \left[1-\frac {(n-2)^2}{n^2} \rho^2 \right]} {1-\frac {n-4}{n} \rho^2}}. \notag
\end{align}
In the same way,
\begin{align}\label{eqn:term3}
III ~=~& \frac {2(n-2)}{n+1} \left[1-\frac {(n-2)^2}{n^2} \rho^2 \right]^{\frac {n+1}{2}} \left(1-\frac {n-2}{n} \rho^2 \right)^{-n} \\
& \quad \times \hyperg {\frac {n}{2}} {\frac {n+1}{2}} {\frac {n+3}{2}}
 {\frac {\rho^2 \left[\frac {(n-2)^2}{n^2} \rho^2 -1\right]} {(1-\frac {n-2}{n} \rho^2 )^2}} \notag\\
=~& \frac {2(n-2)}{n+1} \left[1-\frac {(n-2)^2}{n^2} \rho^2 \right]^{\frac {n+1}{2}} \left(1-\frac {n-4}{n} \rho^2 \right)^{-\frac {n}{2}} \notag\\
& \quad \times  \hyperg {1} {\frac {n}{2}} {\frac {n+3}{2}}
 {\frac {\rho^2 \left[1-\frac {(n-2)^2}{n^2} \rho^2 \right]} {1-\frac {n-4}{n} \rho^2}}.\notag
\end{align}
Applying Gauss' contiguous relation (\cite[p.103, (38)]{EMOT53a})
\begin{equation*}
(c-b)\, z \hyperg {a}{b}{c+1}{z} ~=~ c \hyperg {a-1}{b}{c}{z} - c\, (1-z) \hyperg {a}{b}{c}{z},
\end{equation*}
with $a$=1, $b=\frac {n}{2}$, $c=\frac {n+1}{2}$ and
\[
z=  \frac {\rho^2 \left[1-\frac {(n-2)^2}{n^2} \rho^2 \right]} {1-\frac {n-4}{n} \rho^2},
\]
we have
\begin{align*}
&\hyperg {1}{\frac {n}{2}}{\frac {n+3}{2}}{\frac {\rho^2 \left[1-\frac {(n-2)^2}{n^2} \rho^2 \right]} {1-\frac {n-4}{n} \rho^2}}\\
& \qquad \quad ~=~ \frac {n+1}{\rho^2} \frac {1-\frac {n-4}{n} \rho^2} {1-\frac {(n-2)^2}{n^2} \rho^2 }\\
& \qquad \qquad - \frac {n+1}{\rho^2} \frac {\left(1-\frac {n-2}{n} \rho^2\right)^2} {1-\frac {(n-2)^2}{n^2} \rho^2 }
\hyperg {1}{\frac {n}{2}} {\frac {n+1}{2}}
{\frac {\rho^2 \left[1-\frac {(n-2)^2}{n^2} \rho^2 \right]} {1-\frac {n-4}{n} \rho^2}}.
\end{align*}
Substituting this into \eqref{eqn:term3}, we obtain
\begin{align}\label{eqn:termIII}
III ~=~ & \frac {2(n-2)}{\rho^2} \left[1-\frac {(n-2)^2}{n^2} \rho^2 \right]^{\frac {n-1}{2}}
\left(1-\frac {n-4}{n} \rho^2 \right)^{\frac {2-n}{2}} \\
&\quad  - \frac {2(n-2)}{\rho^2}  \left[1-\frac {(n-2)^2}{n^2} \rho^2 \right]^{\frac {n-1}{2}}
\left(1-\frac {n-4}{n} \rho^2 \right)^{-\frac {n}{2}} \notag\\
& \qquad \quad \times \left(1-\frac {n-2}{n} \rho^2 \right)^2 \hyperg {1} {\frac {n}{2}} {\frac {n+1}{2}}
 {\frac {\rho^2 \left[1-\frac {(n-2)^2}{n^2} \rho^2 \right]} {1-\frac {n-4}{n} \rho^2}}.\notag
\end{align}
Summing up \eqref{eqn:termI}, \eqref{eqn:termII} and \eqref{eqn:termIII} leads to the desired equality \eqref{eqn:2ndderofPhi}.

\section{Proof of Proposition \ref{lem:technical}}

The proof of Proposition \ref{lem:technical} is rather lengthy, we only sketch it. 

Write
\[
\varphi(t):=\frac {t \left[1-\frac {(n-2)^2}{n^2} t\right]} {1-\frac {n-4}{n} t}, \qquad t\in [0,1].
\]
It is easy to check that the inequality \eqref{eqn:technical}  is equivalent to
\begin{align}
\varphi(t)^{\frac {n-1}{2}} \left[1-\varphi(t)\right]^{\frac {1}{2}} \hyperg {1} {\frac {n}{2}} {\frac {n+1}{2}} {\varphi(t)}
~>~ \frac {t^{\frac {n-1}{2}} \left[ 1-\frac {(n-2)^2}{n^2} t \right]^{\frac {n-3}{2}}\left[ 1-\frac {(n-2)(n-3)}{n^2} t \right]}
{\left( 1-\frac {n-4}{n} t \right)^{\frac {n-2}{2}} \left[ 1-\frac {(n-2)(n-3)}{n(n-1)} t \right]}. \notag
\end{align}
So we define
\begin{align}
\Psi(t) ~:=~& \varphi(t)^{\frac {n-1}{2}}  \left[1-\varphi(t)\right]^{\frac {1}{2}} \hyperg {1} {\frac {n}{2}} {\frac {n+1}{2}} {\varphi(t)} \\
& \quad -~ \frac {t^{\frac {n-1}{2}} \left[ 1-\frac {(n-2)^2}{n^2} t \right]^{\frac {n-3}{2}}\left[ 1-\frac {(n-2)(n-3)}{n^2} t \right]}
{\left( 1-\frac {n-4}{n} t \right)^{\frac {n-2}{2}} \left[ 1-\frac {(n-2)(n-3)}{n(n-1)} t \right]}, \qquad t\in [0,1],  \notag
\end{align}
and claim that $\Psi(t) > 0$ for all $t\in (0,1)$. Since $\Psi(0)=0$, it suffices to show that
$\Psi^{\prime} (t)> 0$ for all $t\in (0,1)$.

Applying the elementary relation (\cite[p.102, (23)]{EMOT53a})
\begin{align}
\frac {d}{dz} & \left\{z^{c-a} (1-z)^{a+b-c} \hyperg{a}{b}{c}{z}\right\} \label{eqn:diffhyperg}\\
& \qquad =~ (c-a)\, z^{c-a-1} (1-z)^{a+b-c-1} \hyperg {a-1}{b}{c}{z} \notag
\end{align}
and noting that
\[
\varphi^{\prime}(t) = \frac {\left(1-\frac {n-2}{n} t\right) \left[ 1-\frac {(n-2)(n-4)}{n^2} t \right]}
{\left(1-\frac {n-4}{n} t\right)^2},
\]
we obtain
\begin{align*}
\frac {d}{dt}  & \Bigg\{ \varphi(t)^{\frac {n-1}{2}} \left[1-\varphi(t)\right]^{\frac {1}{2}} \hyperg {1} {\frac {n}{2}} {\frac {n+1}{2}} {\varphi(t)} \Bigg\}\\
&~=~ \frac {n-1}{2}\, \varphi(t)^{\frac {n-3}{2}} \left[1-\varphi(t)\right]^{-\frac {1}{2}} \varphi^{\prime}(t)\\
&~=~ \frac {n-1}{2} t^{\frac {n-3}{2}} \left(1-\frac {n-4}{n} t\right)^{-\frac {n}{2}} \left[1-\frac {(n-2)^2}{n^2} t\right]^{\frac {n-3}{2}}
\left[ 1-\frac {(n-2)(n-4)}{n^2} t \right] \\
& \quad ~=~ t^{\frac {n-3}{2}} \left(1-\frac {n-4}{n} t\right)^{-\frac {n}{2}} \left[1-\frac {(n-2)^2}{n^2} t\right]^{\frac {n-5}{2}}
\left[ 1-\frac {(n-2)(n-3)}{n(n-1)} t \right]^{-2} \notag \\
& \qquad \times \Bigg\{ \frac {n-1}{2} - \frac {12 -46 n + 56 n^2 - 26 n^3 + 4n^4}{2 n^2(n-1)}\, t \notag\\
&\qquad \hspace{1.6cm} + \frac {(n-2)(-8 + 62 n - 119 n^2 + 87 n^3 - 27 n^4 + 3 n^5)}{ n^4 (n-1)}\, t^2 \notag\\
&\qquad \hspace{1.6cm} - \frac {(n-2)^2 (n-3) (32 - 108 n + 98 n^2 - 34 n^3 + 4 n^4)}{2 n^5 (n-1)}\, t^3\notag\\
&\qquad \hspace{1.6cm} + \frac{(n-2)^5 (n-3)^2 (n-4)}{2 n^6 (n-1)}\, t^4 \Bigg\}.
\end{align*}
Also, a direct calculation yields
\begin{align*}
\frac {d}{dt} & \left\{ \frac {t^{\frac {n-1}{2}} \left[ 1-\frac {(n-2)^2}{n^2} t \right]^{\frac {n-3}{2}}\left[ 1-\frac {(n-2)(n-3)}{n^2} t \right]} {\left( 1-\frac {n-4}{n} t \right)^{\frac {n-2}{2}} \left[ 1-\frac {(n-2)(n-3)}{n(n-1)} t \right]} \right\}\\
& \quad ~=~ t^{\frac {n-3}{2}} \left(1-\frac {n-4}{n} t\right)^{-\frac {n}{2}} \left[1-\frac {(n-2)^2}{n^2} t\right]^{\frac {n-5}{2}}
\left[ 1-\frac {(n-2)(n-3)}{n(n-1)} t \right]^{-2} \notag \\
& \qquad \times \Bigg\{ \frac {n-1}{2} - \frac{10 -49 n +57 n^2 -26 n^3 +4 n^4}{2 n^2 (n-1)}\, t\\
&\qquad \hspace{1.6cm} +\frac {(n-2) (-12 + 75 n - 126 n^2 + 88 n^3 - 27 n^4 + 3 n^5)}{ n^4 (n-1)}\, t^2\\
&\qquad \hspace{1.6cm} -\frac {(n-2)^2 (n-3) (48 -116 n +99 n^2 -34 n^3 +4 n^4)}{2 n^5 (n-1)}\, t^3\\
&\qquad \hspace{1.6cm} +\frac {(n-2)^5 (n-3)^2 (n-4)}{2 n^6 (n-1)}\, t^4 \Bigg\}.
\end{align*}
It follows that
\begin{align*}
\Psi^{\prime}(t) =~& \frac {t^{\frac {n-1}{2}}}{2n^5 (n-1)} \left(1-\frac {n-4}{n} t\right)^{-\frac {n}{2}}
\left[1-\frac {(n-2)^2}{n^2} t\right]^{\frac {n-5}{2}} \left[ 1-\frac {(n-2)(n-3)}{n(n-1)} t \right]^{-2} \notag \\
& \quad \times \Big[n^3(n^2-3n-2) - 2n (n-2) (n-4) (n^2 -3 n + 1)\, t \\
&\quad \hspace{3.1cm} +(n-2)^2 (n-3) (n-4)^2\, t^2\Big].
\end{align*}
It is easily seen that, when $n\geq 4$, the quadratic polynomial
\[
n^3(n^2-3n-2) - 2n (n-2) (n-4) (n^2 -3 n + 1)\, t
+(n-2)^2 (n-3) (n-4)^2 \, t^2
\]
is always positive on the interval $[0,1]$. Consequently $\Psi^{\prime}(t)>0$ for all $t\in [0,1]$.
This completes the proof.

\end{document}